\documentclass[aps,prl,showpacs,floatfix,superscriptaddress]{revtex4}%
\usepackage[latin1]{inputenc}
\usepackage{amsmath}
\usepackage{amsfonts}
\usepackage{amssymb}
\usepackage{mathrsfs}
\usepackage{lipsum}
\usepackage[pdftex]{graphicx}
\usepackage{color}

\begin{document}
\title{Additive noise quenches delay-induced oscillations}
\date{\today}
\author{J\'er\'emie Lefebvre}
\email{Jeremie.Lefebvre@unige.ch}
\affiliation{D\'epartement de Neurosciences Fondamentales, University of Geneva, CMU , 1 Rue Michel Servet, 1211 Geneva 4, Switzerland}
\author{Axel Hutt}
\affiliation{INRIA CR Nancy - Grand Est, 615 Rue du Jardin Botanique, 54602 Villers-les-Nancy Cedex, France}
\begin{abstract}
\noindent
Noise has significant impact on nonlinear phenomena. Here we demonstrate that, in opposition to previous assumptions, additive noise interfere with the linear stability of scalar nonlinear systems when these are subject to time delay. We show this by performing a recently designed time-dependent delayed center manifold (DCM) reduction around an Hopf bifurcation in a model of nonlinear negative feedback. Using this, we show that noise intensity must be considered as a bifurcation parameter and thus shifts the threshold at which emerge delay-induced rhythmic solutions.
\end{abstract}
\maketitle
\section{Introduction}
\noindent
Noise has a strong influence on the dynamics of systems subject to time-delayed retroaction. Random fluctuations impact laser dynamics\cite{LangKobayashi}, gene regulatory networks \cite{Bratsunetal2005}, postural and mechanical control \cite{Eurich1996}, which all exhibit retarded self-interactions. In neuroscience, the propagation and time evolution of neural activity has been shown to vary according to combined action of sustained stochastic driving and recurrent feedback, key components of brain circuits \cite{Deco09}.  As such, proper characterization of systems away from the stationnary regime is essential in order to understand how these elements conspire and give rise to the complex phenomena we observe in experiments.\newline
\\
\noindent
Systems evolving under mixtures of noise and nonlinear feedback express divergent properties from their deterministic counterparts, where noise has been shown to provoke or prevent bifurcations \cite{HorsthemkeLefever}. However, such noise-induced stability shifts in scalar systems have are thought to be caused only by a subgroup of stochastic influences, namely multiplicative or parametric noise. This specific type of noise has indeed been shown to be the source of state-transitions in nonlinear systems, while its implication on the stability of delay equations is a well-documented phenomenon \cite{Longtin90,Frank01bpre,Gaudreault2010}.\newline
\\
\noindent
In the past numerous studies have determined or assumed that additive noise, that is stochastic fluctuations independent of the state variable, has no impact on the stability of nonlinear systems \cite{Schumaker}, even in the presence of delays \cite{MackeyNechaeva1995}. This implies that the deterministic expression of the dynamics fully exposes the mean behavior of an otherwise stochastic system. This view has since been challenged in mean-field type systems \cite{Shiino,RodriguezTuckwell1998} and physical pattern forming systems \cite{HuttPhysicaD2008} . More recently, it has been demonstrated that additive noise stabilizes delayed nonlinear systems near non-oscillatory instabilities \cite{HuttEPL2012,LefebvreCHAOS2012}. In this Letter, we show that additive stochastic fluctuations cause a shift in the eigenspectrum in a nonlinear scalar delayed differential equation(DDE) in the vicinity of an Hopf bifurcation and de facto shapes delay-induced rhythms. By doing so, we state that despite the current belief, additive noise changes significantly the stability of delayed oscillatory nonlinear systems and that the additive noise intensity has to be considered as a bifurcation parameter.  \newline
\\
\noindent
To demonstrate this, we adopt the following strategy. We first derive a time-dependent delayed center manifold (DCM) representation of a generic scalar DDE unfolded around a supercritical Hopf bifurcation. Then, adiabatic elimination is performed based on the characteristic separation of time scales involved in the vicinity of a non-hyperbolic fixed point. We show that the explicit time-dependence of the DCM reveals noise-dependent components reponsible for the depletion of rhythmic solutions for negative nonlinearities, and consequently exposes the impact of additive noise on the system's equilibrium and its stability.  We then apply our theory to reconstruct the associated stochastic bifurcation diagram of the ensemble averaged noisy solutions, and validate it via numerical simulations.\newline

\section{Model}
\noindent
The work is motivated by non-homogeneous retarded equations with the generic form 
\begin{equation}
\dot{x}(t)=-x(t)+\mathcal \gamma f(x(t-\tau))+I(t),
\label{DDE0}
\end{equation}
\noindent
exhibiting a smooth nonlinearity $f(x)$ and some additive random fluctuation $ I(t)$, in the Itô sense.  Let us consider the expansion of Eq.\ (\ref{DDE0}) around the stationary state $x_o=\gamma f(x_o)$ up to third order, which yields the polynomial system 
\begin{equation}
\dot{u}(t)= -u+\eta u_{\tau}+ \kappa u_{\tau}^2+ \nu u_{\tau}^3 + \sqrt{2D}\xi(t)~,
\label{DDE}
\end{equation}
\noindent
where $\eta=\gamma f'(x_o)$, $\kappa=\gamma f''(x_o)/2!$ and $\nu=\gamma f'''(x_o)/3!$ are non-zero constants, with the small deviations from the stationary state  $u=u(t)$ and $u_{\tau}=u(t-\tau)$. Here, the additive random uncorrelated fluctuations are defined by $I(t)=\sqrt{2D}\xi(t)$ with the noise intensity $D>0$ such that $\langle\xi(t')\xi(t)\rangle=\delta(t-t')$ where $\langle\cdot\rangle$ denotes the 
ensemble average.\newline
\\
\noindent
Figure \ref{effect} shows significant differences between the dynamics of Eq.(\ref{DDE}) for $D=0$ and $D>0$ close to the bifurcation point. While virtually masked by local fluctuations in single trial trajectories, the investigation of the ensemble average reveals that additive noise stabilizes the oscillations by preventing the bifurcation. This is caused by a shift of the eigenspectrum performed by noise and delay. The impact of noise decreases as $|\epsilon|$ increases i.e. away from the bifurcation, the time-scale separation vanishes and the system's ensemble average becomes somewhat passive to addive noise.\newline

\begin{widetext}
\begin{center}
\begin{figure}[h!]
\includegraphics[scale=0.08]{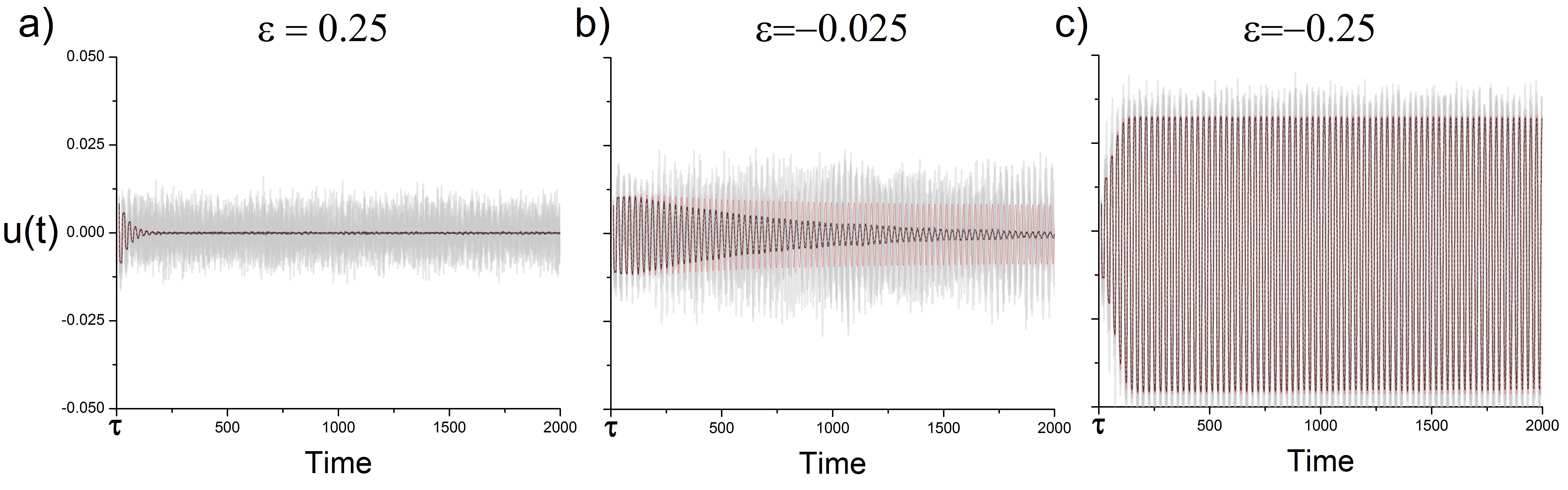} 
\caption{\textbf{Impact of additive noise on the stability of delay induced oscillations}. From panel a) to c), the system (\ref{DDE}) undergoes a supercritical Hopf bifurcation as the control parameter $\epsilon=\eta-\eta_c$ is varied from $0.25$ to $-0.25$ and crosses the Hopf threshold located at $\epsilon=0$. Slightly above the bifurcation in panel b), one observes deviations between the noise-free solution $u(t)$(red curve; $D=0$) and the ensemble averaged solution $\langle u(t)\rangle$ (dark grey curve) over independent realizations of the noise. While the noise-free system exhibits stable oscillations of fixed amplitudes, noise prevents the bifurcation and mean solution's amplitude decreases (panel b)). The deviation between noisy and noise-free dynamics vanishes as one moves away below (panel a)) or above (panel c)) the bifurcation threshold. This effect is not easily observable amongst the single-trial realizations (light gray). DCM corrected solutions capture the effect of the additive noise on stability and reproduces the dynamics of the ensemble average (dashed curve, almost undistinguishable from the ensemble average).  Parameters are $D=0.00001$, $\gamma=-0.05$, $\tau=12.0$, $\eta_c=-1.03$ and $f(x)=(1+\text{exp}[-60\cdot x])^{-1}$. The time integration followed a Euler-Maruyama scheme with $10^4$ iterations and time step $dt=0.1$.  \label{effect} }
\end{figure}
\end{center}
\end{widetext}

\section{Unfolding around the Hopf bifurcation}
\noindent
In order to explain this novel phenomenon, the subsequent paragraphs at first neglect noise and follow the standard analysis of DCM reduction for deterministic delay equations \cite{Hale93book,Magalhaes,PelsterSIAM,Campbell}. Analyzing the stationary state's linear stability, we consider just the linear terms in Eq.~(\ref{DDE}) 
\begin{eqnarray}
\dot{u}(t)= -u+\eta u_{\tau}\label{DDE_lin}
\end{eqnarray}
and focus on the Hopf instability assuming $u(t)=u_0\exp(\lambda t),~\lambda\in{\mathcal C},~u_0=$const. Then there is a critical $\lambda_c=\pm iw_c$ that satisfies
\begin{equation}
\lambda_c=-1+\eta_c e^{\lambda_c \tau_c}~,
\label{eigen1}
\end{equation}
\noindent 
where $w_c=\sqrt{\eta_c^2-1}$ is the corresponding Hopf frequency. Introducing the control parameter $\epsilon = \eta-\eta_c$, it is convenient to re-write Eq.~(\ref{DDE}) as 

\begin{equation}
\dot{u}(t)= L(u,\eta_c)+F(\epsilon,u)+ \sqrt{2D}\xi(t),\ \dot{\epsilon}(t)=0\ ,
\label{DDEAUG}
\end{equation}
with $L(u,\eta_c)\equiv-u+\eta_c u_{\tau}$ and $F(\epsilon,u)\equiv\epsilon u_{\tau}+ \kappa u_{\tau}^2+ \nu u_{\tau}^3$. In this formulation, the Hopf bifurcation occurs for $\epsilon=0$.\newline

\section{DCM reduction neglecting noise}
\noindent
Equation (\ref{DDE_lin}) and the characteristic equation (\ref{eigen1}) define eigenvectors $\Phi(\theta)=(\phi_1(\theta),\phi_2(\theta))$ where $\theta$ is a parametrization of the delay i.e. $-\tau\leq\theta\leq 0$. They span the center subspace and are determined by $\pm iw_c \phi_j(0) = L(\phi_j,\eta_c)=-\phi_j(0)+\eta_c \phi_j(-\tau)|j=(1,2)$ yielding 
\begin{equation}
\Phi(\theta)=(\phi_1(\theta),\phi_2(\theta))=(e^{i w_c \theta},e^{-i w_c \theta})~.
\end{equation}
\noindent 
The associated adjoint basis of the center subspace $\Psi(s)=(\psi_1(s),\psi_2(s))^\intercal = (d e^{-i w_c s}, \bar{d} e^{i w_c s})^\intercal $, $d\in \mathbb C$, obeys $(\Phi,\Psi)=\mathbb I$, where $(a(\theta),b(\theta))$ is a bilinear form defined by [Lunel]
\begin{equation}
\textbf{(}a,b\textbf{) }\equiv a(0)b(0) - 
\int_{-\tau}^{0}\int_{0}^{\theta}a(\xi-\theta)b(\xi)d\xi[d\alpha(\theta)]~.
\label{bilinearform}
\end{equation}
Here, $[d\alpha(\theta)]=(-\delta(\theta)+\eta_c\delta(\theta+\tau))d\theta$.  According to the center manifold theorem, phase space can be split into center ($C$) and stable ($S$) parts i.e. $C\oplus S$, such that the stables modes converge to the center manifold. We thus have 
\begin{eqnarray*}
u_t(\theta)=\Phi(\theta)\textbf{z}+h(\textbf{z},\theta)\in C
\end{eqnarray*} 
for some stable manifold $h$, where $\textbf{z}=[z_1(t),z_2(t)]$ is the amplitude on the center manifold. The dynamics of Eq.\ (\ref{DDEAUG}) projected onto the linear center subspace is consequenly given by the order parameter equation
\begin{equation}
\dot{\textbf{z}}(t) = B\textbf{z}+\Psi(0) F[\epsilon,\Phi\textbf{z}+h],\ \dot{\epsilon}(t)=0~,
\label{OPE1}
\end{equation}
\noindent
where the diagonal matrix B is defined by
\[B= \left( \begin{array}{ccc}
i w_c & 0 \\
0 & -i w_c  \end{array} \right).\]

\section{Stochastic DCM reduction}
\noindent
Now we extend this standard DCM analysis close to a Hopf instability by introducing additive noise i.e. $D>0$. We treat the additive noise similar to the nonlinearity which renders the stable and center manifolds time-dependent. Hence Eq.\ (\ref{OPE1}) reduces to
\begin{eqnarray}
\dot{\textbf{z}}(t) &=& B\textbf{z}+\Psi(0) F[\epsilon,\Phi(\theta)\textbf{z}+h(\textbf{z},\theta,t)]\nonumber\\
&&+\sqrt{2D}\Psi(0)\xi(t),\ \dot{\epsilon}(t)=0 ~.\label{OPE}
\end{eqnarray}
\noindent
The separable ansatz for the \textit{time dependent} DCM close to a non-oscillatory instability \cite{HuttEPL2012,Boxler} reads
\begin{equation}
h(\textbf{z},\theta,t) = h(\textbf{z},\theta)+h_t(\theta,t)\ .
\label{h}
\end{equation}
\noindent
Here $h_t(\theta,t) = (1-P_{c})H(t+\theta)$ is independent of $z$ and evolves on the linear stable subspace. The operator $P_{c}$ projects onto the center subspace spanned by $\Phi(\theta)$\cite{LefebvreCHAOS2012} and $H(t)\in{\mathcal C}$ obeys 
\begin{equation}
\dot{H(t)} = L(H,\eta_c)+\sqrt{2D}\xi(t).
\label{H}
\end{equation}
\noindent
The separation ansatz in Eq.\ (\ref{h}) provides a good approximitive description of delayed noisy dynamics near instabilities in regimes where $\Phi(\theta)\textbf{z}+h(\textbf{z})$ evolves slowly with small amplitude. In the vivinity of an oscillatory instability, this corresponds to regimes where the emerging oscillation has a small frequency and small amplitude, which occurs when $\tau$ is taken large enough. This is the regime we investigate here. This choice is further in line with past numerical findings where the impact of noise was shown to be inversely proportional to the delay size \cite{Longtin91}. \newline
\noindent
\\ 
The state probability density $p(\textbf{z},h_t)$ of the system evolving according to Eq.\ (\ref{OPE}) considers $h_t$ as a time-dependent random variable and obeys the Fokker-Planck equation,
\begin{eqnarray}
&&\frac{\partial p(\textbf{z},h_t)}{\partial t}= 
\psi_1(0)\left\lbrace D\frac{\partial p^2(\textbf{z},h_t)}{\partial^2 z_1}\right.\nonumber\\
&&\left.-\frac{\partial}{\partial z_1}\left(i w_c z_1+F[g(\textbf{z})+h_t]\right)p(\textbf{z},h_t)\right \rbrace +\psi_2(0)\bigg\lbrace  \nonumber\\
&&\left.D\frac{\partial p^2(\textbf{z},h_t)}{\partial^2 z_2}-\frac{\partial}{\partial z_2}\left(-i w_c z_2+F[g(\textbf{z})+h_t]\right)p(\textbf{z},h_t) \right \rbrace\nonumber\\ 
&&\label{FPE}
\end{eqnarray}
\noindent
where $g(\textbf{z})=\Phi(\theta) \textbf{z}+h(\textbf{z},\theta)$ is the time-independent contribution to the center manifold. Here the dependence in $\theta$ and the trivial dynamics of the unfolding parameter $\epsilon$ have been dropped for clarity. \newline
\\
\noindent
The characteristic time scale separation involved near non-hyperbolic fixed points implies that the time-dependent functional $h_t(t)$ and the noisy input $\xi(t)$ are much smaller than the center modes $\textbf{z}$ and retain the time-scale separation of the unperturbed system. By virtue of the separable ansatz in Eq.\ (\ref{h}) for which the temporal correction $h_t$ is independent of $\textbf{z}$, the adiabatic elimination \cite{DroletVinals2001,HuttPhysicaD2008} of the fast time dependent terms applies by performing an ensemble average of $p(\textbf{z},h_t)$ over $h_t$ for all $\theta$ in Eq.~(\ref{FPE}) using the associated density $q(h_t)$. Specifically, for weak noise one separates the slow and fast evolution by   $p(\textbf{z},h_t)=p_0(\textbf{z},t)q(h_t)$ and expand the nonlinearity $F$ with respect to $h_t$ around the center modes. 
Assuming stationarity of $h_t$ and a Gaussian profile for $p(h_t)$ with vanishing mean and variance $\sigma^2 = \sigma^2(D) = \int_{0}^{T} h^2_t(\theta,s)p(h_t(\theta,s))ds>0 $ for a short time window $T\ll 2\pi/w_c$ shorter than the time scale of the center dynamics
\begin{eqnarray*}
&&\int_{t}^{t+T} F[g(\textbf{z}(t))+h_t(\theta,s)]p(h_t(\theta,s))ds \approx\\
&& \epsilon g_\tau(\textbf{z}(t))+ \kappa \left(g_\tau^2(\textbf{z}(t))+ \sigma^2\right)
+\nu \left(g_\tau^3(\textbf{z}(t))+3 g_\tau(\textbf{z}(t))\sigma^2\right).
\end{eqnarray*}
with $g_\tau(\textbf{z}(t))=\Phi(-\tau)\textbf{z}(t)+h(\textbf{z}(t),-\tau)$. 
The resulting adiabatic Fokker Planck equation determines the temporal evolution of the probability density $p_0(\textbf{z},t)$ whose mean value $\bar{\textbf{z}}$ obeys the non-autonomous order parameter equation
\begin{eqnarray}
&&\dot{\bar{\textbf{z}}}(t) = c_o+B\bar{\textbf{z}}(t)+\Psi(0) F(\epsilon^*,\Phi(\theta)\bar{\textbf{z}}(t)+h(\bar{\textbf{z}}(t),\theta)),\nonumber\\
&&\label{OPE_cor}
\end{eqnarray}
with the effective control parameter $\epsilon^*(D)=\epsilon+\mu(D)$, $c_o(D)=\kappa \sigma^2(D)$ and $\mu(D)=3 \nu \sigma^2(D) $.

\noindent 
A comparison between the full and reduced description reveals the mechanism by which additive noise alters the properties of the noise-free system. To see this, note that Eq.\ (\ref{OPE_cor}) may be viewed as the reduced dynamics of the non-autonomous DDE
\begin{equation}
\label{DDEABS}
\dot{Z}(t) = c_o(D)+L(Z(t),\eta_c)+F(\epsilon^*(D),Z(t))
\end{equation}
\noindent
This equation is our major result. By direct comparison with Eq.\ (\ref{DDEAUG}), two effects can be seen in this system with cubic nonlinearity:  \textit{i}) Since $c_o\neq 0$, for $D>0$ the system has a new equilibrium repositioned in phase space and, \textit{ii}) the noise shifts the control parameter by a factor of $\mu$ and thus alters the stability of the system. Note that Eq.\ (\ref{DDEABS}) converges to Eq.\ (\ref{DDEAUG}) when $D=0$. Summarizing, the sign and manitude of the nonlinear parameters $\kappa$ and $\nu$ in the expanded system fully determine the qualitative impact of the additive perturbations on stability. This finding demonstrates the purely nonlinear nature of phenomenon shown in Fig.~\ref{effect}.

\noindent
To test our theory, we have numerically investigated the behavior of the DDE in Eq.(\ref{DDE}) with $D=0$ and $D>0$ where the delayed nonlinearity is taken to be negative ($\gamma<0$). 
To compare the original dynamics with the DCM corrected system in Eq. (\ref{DDEABS}),  
we have chosen a regime of large delays where the center modes oscillate at a small frequency and where the impact of noise on the dynamics is known to be weak \cite{Longtin91}. To determine $\sigma^2$, we have integrated numerically the linear system in Eq. (\ref{H}), applied the corresponding definition of $h_t$, and computed the stationnary variance $\sigma^2 \approx (1-2 \text{Re}(d)\cos(w_c \tau))^2\langle H^2\rangle$, which holds true for $D>0$ small. The constant $d$ is the normalization constant of the adjoint basis $\Psi(s)$.

\begin{center}
\begin{figure}[h!]
\includegraphics[width=8cm]{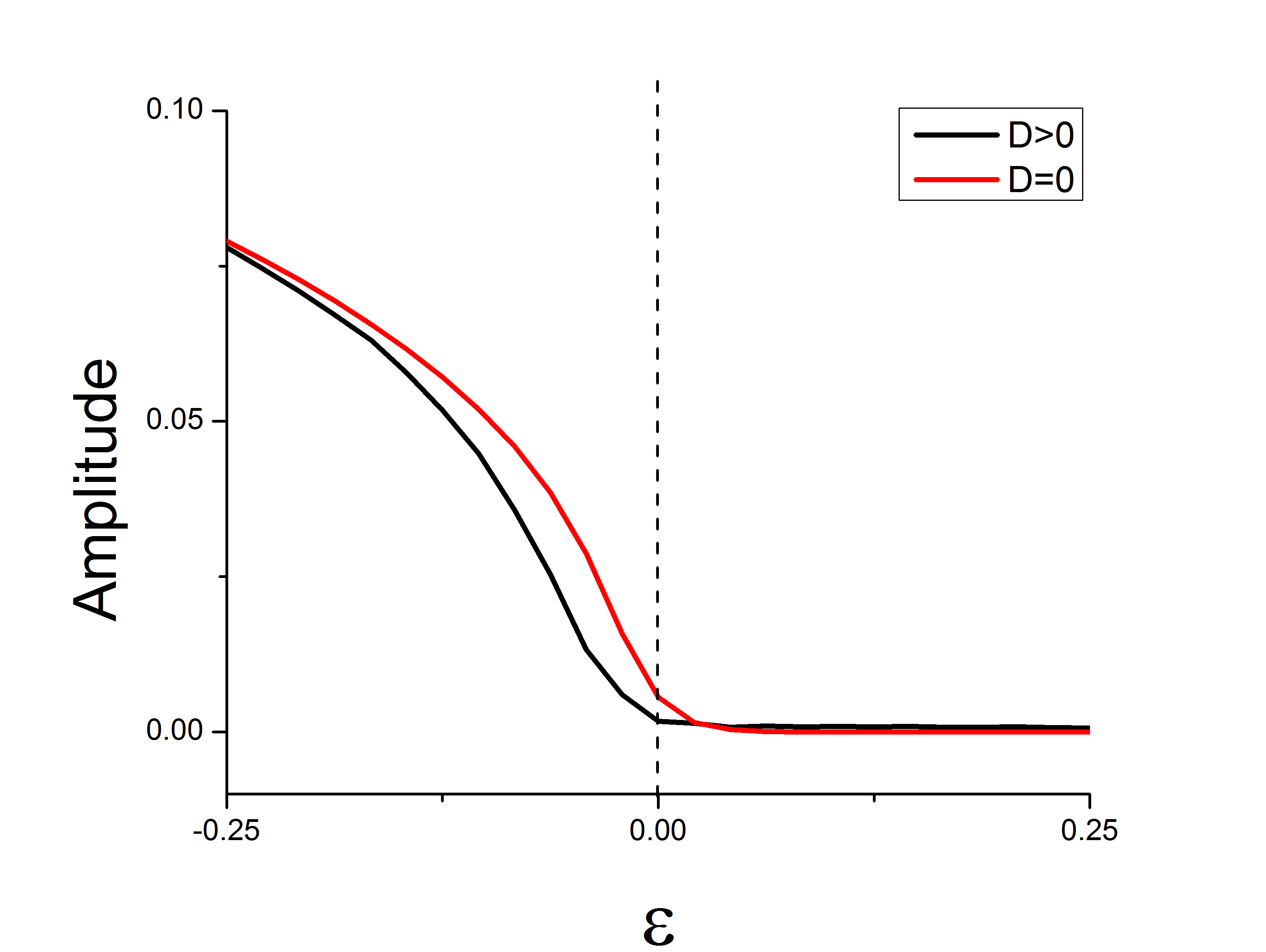} 
\caption{\textbf{Bifurcation diagram for oscillation amplitudes subjected to the control parameter.} The amplitudes of the oscillatory solutions are computed by numerical integration of the original stochastic DDE in Eq.\ (\ref{DDE}) for $D=0$(red) and $D=0.00001$(black). We observe that the two curves are well distinguishable in the neighborhood of the bifurcation point at $\epsilon=0$, but coalesce otherwise. Other parameters are taken from Fig.\ \ref{effect}. \label{shift} }
\end{figure}
\end{center}

\noindent
While set in  such a regime of slow oscillations, Fig.~\ref{bifplot} shows that additive noise shifts the Hopf bifurcation to smaller values and thus stabilizes the system. The corrected system (\ref{DDEABS}) captures very well the dynamics altered by additive noise and further allows to predict the position of the noise-shifted Hopf bifurcation in parameter space. Indeed, the DCM reduction scheme captures both qualitatively and quantitatively the effect of additive noise on the dynamics of the original equation in the neighborhood of the bifurcation whenever $D$ and $|\epsilon|$ are chosen small. It further reveals the characteristic quenching action of the random fluctuation on the ensemble averaged solutions experienced in presence negative delayed feedback. 
\begin{widetext}
\begin{center}
\begin{figure}[h!]
\includegraphics[scale=0.08]{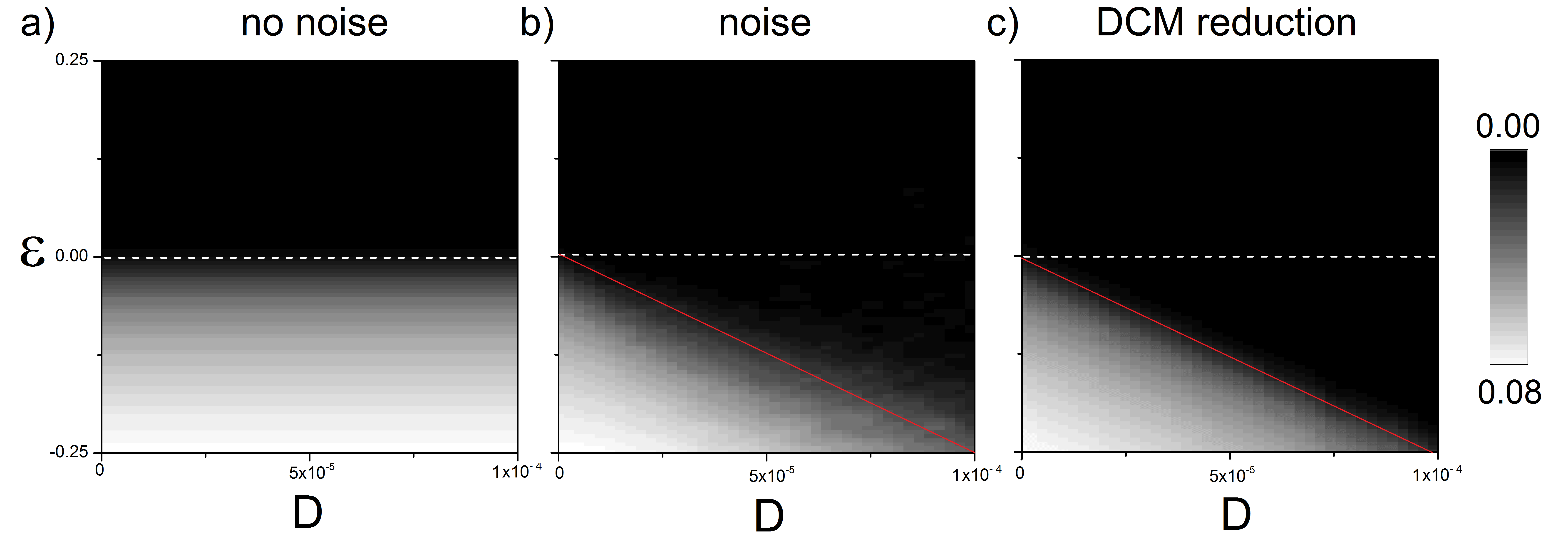} 
\caption{\textbf{Delayed stochastic Hopf bifurcation diagram}. Amplitudes of the nascent oscillatory solutions as a function of noise intensity ($D$) and distance from the bifurcation threshold $\epsilon$. The Hopf threshold is located at $\epsilon=0$ (dashed line across plots for comparison). a) Generic bifucation diagram in absence of noise ($D=0$). b) As noise is injected additively in the system, the onset of stable oscillatory solutions retreats in parameter space compared to the noise-free case depicted in panel a), giving rise to significantly different behavior. The panel shows the amplitudes of the ensemble averaged solutions $\langle u(t)\rangle$ of Eq. (\ref{DDE}). c) Amplitudes of the corrected system (\ref{DDEABS}) capture the additive noise effect. Deviations between the noisy diagram and the DCM reduction prediction (between panels b) and c)) appear as $D$ and $|\epsilon|$ increase pointing to the approximative nature of the separation ansatz (\ref{h}). Other parameters are taken from Fig.~\ref{effect}. \label{bifplot} }
\end{figure}
\end{center}
\end{widetext}
\noindent
One may ask whether the effect revealed, i.e. the shift of the control parameter and the oscillation amplitude, is prominent enough to play a significant role in system's dynamics. Figure~\ref{shift} shows oscillation amplitudes subjected to the control parameter $\epsilon$ for two noise levels. It reveals that the oscillation amplitudes are practically indistinguishable for $\epsilon<-0.25$ or $\epsilon>0.02$, i.e. away from the original deterministic bifurcation threshold. Hence the additive noise effect on the ensemble average of the system activity is observable near deterministic stability threshold but disappears far from the bifurcation point. 

\section{Discussion}
\noindent
Our results reveal that additive noise perturbs the stability of scalar delayed systems and shifts the instability threshold in parameter space, demonstrating that noise intensity must be considered as a bifurcation parameter. Using time-dependent DCM reduction, we have exploited the characteristic time-scale separation emerging near non-hyperbolic fixed points to perform adiabatic elimination. For the cubic system considered, our approach reveals the effect of additive noiseon linear and cubic terms. This perturbs the stability of the systems and shifts the onset of oscillatory solutions.\newline
\\
\noindent
This novel finding has several implications. First, the effect suggests a general mechanism for on-line suppression of rhythmic activity by external stochastic driving in presence of inhibitory feedback, despite purely deterministic model parameters. We conjecture that this mechanism might provide support to experimental paradigms seeking to supress undesired or pathological rhythmic regimes in recurrent systems. \newline
\\
\noindent
The work considers the supercritical Hopf bifurcation for a delayed nonlinearity with negative gain ($\gamma<0$). Although our theory is not specific to this case, we note that in the vicinity of a subcritical instability, additive noise causes the solutions to diverge and that, in most practical applications, the supercritical case prevails. This case nonetheless suggests that noise in such systems might have the opposite effect: it might enhance oscillatory behavior. This remains to be shown in future work.

\nocite{*}
\bibliographystyle{unsrt}

\end{document}